\documentclass[12pt]{amsart}

\begin{document}
\numberwithin{equation}{section}

\def\1#1{\overline{#1}}
\def\2#1{\widetilde{#1}}
\def\3#1{\widehat{#1}}
\def\4#1{\mathbb{#1}}
\def\5#1{\frak{#1}}
\def\6#1{{\mathcal{#1}}}

\def\C{{\4C}}
\def\R{{\4R}}
\def\N{{\4N}}
\def\Z{{\4Z}}

\def\bC{{\4C}}
\def\bR{{\4R}}
\def\bN{{\4N}}
\def\bZ{{\4Z}}
\def \bP{{\4P}}

\title[Invariants of CR Manifolds in $\bC^2$]{Invariants and Umbilical Points on Three Dimensional CR Manifolds embedded in $\bC^2$}
\author[P. Ebenfelt and D. Zaitsev]{Peter Ebenfelt and Dmitri Zaitsev}
\address{D. Zaitsev: School of Mathematics, Trinity College Dublin, Dublin 2, Ireland}
\email{zaitsev@maths.tcd.ie}
\address{Department of Mathematics, University of California at San Diego, La Jolla, CA 92093-0112}
\email{pebenfel@math.ucsd.edu}

\abstract We introduce a new sequence of CR invariant determinants on a three dimensional CR manifold $M$  embedded in $\mathbb C^2$. The lowest order invariant $\det A_3$ represents E. Cartan's 6th order invariant (the umbilical "tensor"), whose zero locus yields the set of umbilical points on $M$. As an application of this new presentation of the umbilical invariant, we show that generic, almost circular perturbations of the unit sphere always contain curves or surfaces of umbilical points. 
\endabstract

\thanks{The first author was supported in part by the NSF grant DMS-1301282.}


\maketitle

\def\Label#1{\label{#1}}


\def\cn{{\C^n}}
\def\cnn{{\C^{n'}}}
\def\ocn{\2{\C^n}}
\def\ocnn{\2{\C^{n'}}}


\def\dist{{\rm dist}}
\def\const{{\rm const}}
\def\rk{{\rm rank\,}}
\def\id{{\sf id}}
\def\aut{{\sf aut}}
\def\Aut{{\sf Aut}}
\def\CR{{\rm CR}}
\def\GL{{\sf GL}}
\def\Re{{\sf Re}\,}
\def\Im{{\sf Im}\,}
\def\span{\text{\rm span}}

\def\norm{|\!|}

\def\codim{{\rm codim}}
\def\crd{\dim_{{\rm CR}}}
\def\crc{{\rm codim_{CR}}}

\def\phi{\varphi}
\def\eps{\varepsilon}
\def\d{\partial}
\def\a{\alpha}
\def\b{\beta}
\def\g{\gamma}
\def\G{\Gamma}
\def\D{\Delta}
\def\Om{\Omega}
\def\k{\kappa}
\def\l{\lambda}
\def\L{\Lambda}
\def\z{{\bar z}}
\def\w{{\bar w}}
\def\Z{{\1Z}}
\def\t{\tau}
\def\th{\theta}

\def\mU{\mathcal U}
\def\AC{\mathcal A\mathcal C}
\def\RP{\mathcal P^{\bR}}
\def\FS{\mathcal F\mathcal S}
\def\FL{\text{FL}}

\emergencystretch15pt
\frenchspacing

\newtheorem{Thm}{Theorem}[section]
\newtheorem{Cor}[Thm]{Corollary}
\newtheorem{Pro}[Thm]{Proposition}
\newtheorem{Lem}[Thm]{Lemma}

\theoremstyle{definition}\newtheorem{Def}[Thm]{Definition}

\theoremstyle{remark}
\newtheorem{Rem}[Thm]{Remark}
\newtheorem{Exa}[Thm]{Example}
\newtheorem{Exs}[Thm]{Examples}
\newtheorem{Que}[Thm]{Question}

\def\bl{\begin{Lem}}
\def\el{\end{Lem}}
\def\bp{\begin{Pro}}
\def\ep{\end{Pro}}
\def\bt{\begin{Thm}}
\def\et{\end{Thm}}
\def\bc{\begin{Cor}}
\def\ec{\end{Cor}}
\def\bd{\begin{Def}}
\def\ed{\end{Def}}
\def\br{\begin{Rem}}
\def\er{\end{Rem}}
\def\be{\begin{Exa}}
\def\ee{\end{Exa}}
\def\bpf{\begin{proof}}
\def\epf{\end{proof}}
\def\ben{\begin{enumerate}}
\def\een{\end{enumerate}}
\def\beq{\begin{equation}}
\def\eeq{\end{equation}}

\section{Introduction}

The motivation behind this paper is a basic question, due to S.-S. Chern and J. K. Moser \cite{CM74}, which roughly asks if there are compact, strictly pseudoconvex CR manifolds of dimension three that do not possess (CR) umbilical points. The answer to the question, in this generality, is well known to be 'yes'. E. Cartan had much earlier \cite{Cartan33} discovered a 1-parameter family of real hypersurfaces $\mu_\alpha\subset \bP^2$ that are compact, strictly pseudoconvex, homogeneous and non-spherical (hence without umbilical points).  A family of 2:1-covers of these CR manifolds are embeddable in $\bC^3$. The universal cover of $\mu_\alpha$ is the sphere $S^3$, and by pulling back the CR structures of $\mu_\alpha$ one obtains a family of CR structures on $S^3$ that do not possess any umbilical points. The latter structures, however, are well known to not be embeddable in $\bC^n$ for any $n$. (See \cite{Isaev06}, \cite{Rossi65}.) A precise question that remains to be answered is then the following:

\begin{Que}\Label{CMQ} Does there exist a compact, strictly pseudoconvex CR manifold $M$ in $\bC^2$ that does not have any umbilical points?
\end{Que}

We point out the analogy with the classical Caratheodory Conjecture (see, e.g., \cite{Hamburger40}, \cite{Ivanov02}) regarding umbilical points on compact surfacces embedded in $\bR^3$, and refer the interested reader to the paper \cite{EDumb15} for a closer discussion of the analogy between the notion of (CR) umbilical points in CR geometry and that of umbilical points in the classicial geometry of surfaces in $\bR^3$.

Let $M=M^{2n+1}$ be a real hypersurface in $\bC^{n+1}$ and $p\in M$. The lowest order (local) invariant of $M$ is its Levi form at $p$, which roughly speaking is a Hermitian form on the tangent space $T^{1,0}_pM$ of $(1,0)$-vectors tangent to $M$ at $p$. A hypersurface is said to be Levi nondegenerate at $p$ if the Levi form is nondegenerate at $p$, and strictly pseudoconvex if the Levi form at $p$ is definite. Levi nondegeneracy of $M$ at $p$ can be detected by the condition $J_p\neq 0$, where $J=J(\rho)$ is
Fefferman's Monge-Ampere operator \cite{Fefferman76}
\beq\Label{FeffJ}
J(\rho):=(-1)^{n+1} \det
\begin{pmatrix}
\rho & \rho_{\bar Z}\\
\rho_{Z} &\rho_{Z\bar Z}
\end{pmatrix}
\eeq
applied to a local defining function $\rho$ of $M$ near $p$; i.e., $M$ is locally given by $\rho=0$, and $d\rho\neq 0$ on $M$. In \eqref{FeffJ}, the notation used is $\rho_{\bar Z}:=(\rho_{\bar Z_1},\ldots,\rho_{\bar Z_{n+1}})$, $\rho_{Z}$ is its conjugate transpose $\rho_Z=(\rho_{\bar Z})^*$, and $\rho_{Z\bar Z}$ is the $(n+1)\times(n+1)$ matrix $(\rho_{Z_k\bar Z_j})$. One of the main contributions in this paper is the introduction of a sequence of invariant determinants, for hypersurfaces in $\bC^2$, that can be viewed as higher order analogs of Fefferman's operator $J$ in \eqref{FeffJ}.

Chern and Moser showed \cite{CM74} that if $M$ is strictly pseudoconvex at $p$, then that there are formal holomorphic coordinates $Z=(z,w)=(z_1,\ldots,z_n,w)$ (convergent if $M$ is real-analytic), vanishing at $p$, such that $M$ can be expressed in Chern--Moser normal form. Rather than describing this normal form precisely here, we simply note that in these coordinates $M$ is expressed as a graph \beq
\Im w=\phi(z,\bar z,\Re w),
\eeq
where the graphing function has the form
\beq\Label{CMnormform}
\phi(z,\bar z,\Re w)=\sum_{j=1}^n|z_j|^2+R_{m}(z,\bar z)+O(m+1).
\eeq
Here, $R_{m}(z,\bar z)$ is a homogeneous Hermitian polynomial of degree $m$ with $m=6$ for $n=1$ and $m=4$ for $n\ge 2$, $O(k)$ signify terms of weight $\geq k$ in $(z,\bar z, \Re w)$, where $(z,\bar z)$ are assigned weight one and $\Re w$ weight two. If $n\geq 2$, then $m=4$ and $R_4(z,\bar z)$ is of bidegree $(2,2)$; $R_4(z,\bar z)$ represents the CR curvature tensor $S_{\alpha\bar\beta\nu\bar\mu}$ at $p=(0,0)$ as its sectional curvature
\beq\Label{R_4}
R_4(z,\bar z)=\sum_{\alpha,\beta,\nu,\mu}S_{\alpha\bar\beta\nu\bar\mu}z_\alpha z_\nu\overline{z_\beta z_\mu}.
\eeq
The CR curvature tensor $S_{\alpha\bar\beta\nu\bar\mu}$ has Hermitian curvature symmetries and zero trace. (The latter is equivalent to $R_4(z,\bar z)$ being a harmonic function of $z$.) The real dimension of the space of such curvature tensors (of Weyl--Bochner type) is $n^2(n-1)(n+3)/3$ (see \cite{sitaramayya73}). From this it follows that the condition of being umbilical at a point $p$ is equivalent to $n^2(n-1)(n+3)/3$ independent real equations on the 4-jet of the CR manifold $M^{2n+1}$ at $p$. For $n=2$, this means 5 equations on a 5-dimensional manifold, but for $n\geq 3$, this is an "overdetermined" system. More precisely, an application of Thom's Transversality Theorem \cite{GG86} (see \cite{BRZ07} for similar arguments) shows that a {\it generic} strictly pseudoconvex hypersurface $M^{2n+1}\subset \bC^{n+1}$ (i.e., in a dense open subset in the compact-open topology of $C^\infty$-mappings $M^{2n+1}\hookrightarrow\bC^{n+1}$) has no umbilical points when $n\geq 3$, and at most isolated umbilical points when $n=2$. Moreover, to further support the statement that umbilical points are rare when $n\geq2$, we mention that Webster has shown \cite{Webster00} that every non-spherical real ellipsoid in $\bC^{n+1}$ is free of umbilical points when $n\geq 2$.


In this paper, we shall consider the case $n=1$, i.e., strictly pseudoconvex hypersurfaces $M=M^3$ in $\bC^2$. In this case, the CR curvature vanishes identically, and the lowest order nontrivial invariant in the Chern--Moser normal form occurs in weight $m=6$, and $R_6(z,\bar z)$ in \eqref{CMnormform} has the form
\beq\Label{R_6}
R_6(z,\bar z)=c_{24}z^2\bar z^4+c_{42}z^4\bar z^4,\quad c_{42}=\overline{c_{24}}\in \bC,
\eeq
where $c_{24}$ represents E. Cartan's "6th order tensor" $Q=Q^1{}_{\bar 1}$ at $p=(0,0)$. The hypersurface $M^3$ is umbilical at $p=(0,0)$ if $c_{24}=0$. Thus, the condition of being umbilical on $M^3$ amounts to two independent real equations in the 6-jet space of CR manifolds in $\bC^2$. Thus, we should expect umbilical points (if they exist!) to form real curves in $M^3$. If we consider the condition that $M^3$ is umbilical at $p$ and the rank of the differential of the two real equations $c_{24}=0$ is $\leq 1$, which amounts to the vanishing of two real determinants, then we obtain an algebraic subvariety in the 7-jet space of codimension $4$. An application of Thom's Transversality Theorem, as above, yields that a {\it generic} strictly pseudoconvex hypersurface $M^3\subset \bC^2$ either has no umbilical points at all or a set of umbilical points that consists of smooth real curves.

\subsection{Summary of Main Results.} The major problem addressed in this paper is that of existence of umbilical points. As is illustrated by the family of examples $\mu_\alpha\subset\bP^2$ and their covers, there are compact three dimensional CR manifolds without umbilical points. We shall focus here on Question \ref{CMQ} described above. Not much is known in general about this problem. It was shown by X. Huang and S. Ji \cite{HuangJi07} that every real ellipsoid in $\bC^2$ must have umbilical points. More recently, it was proved by the first author and S. Duong \cite{EDumb15} that every {\it circular} $M^3\subset \bC^2$ has umbilical points; in fact, it was proved in \cite{EDumb15} that every compact three dimensional CR manifold with a transverse free CR $U(1)$ (circle) action must have umbilical points provided that the Riemann surface $M/U(1)$ has genus $g\neq 1$. As a first step towards answering Question \ref{CMQ} more generally, we shall consider small perturbations $M_\eps\subset \bC^2$ of the unit sphere $M_0=S^3\subset \bC^2$. We shall show that generic {\it almost circular} perturbations $M_\eps$ must possess umbilical points. This is Theorem \ref{ACpertThm}. Precise statements and definitions are given in Section \ref{PertSec}. We also note in Remark \ref{ACrem} that real ellipsoids are not generic in the sense of Theorem \ref{ACpertThm}, but we give separately a new proof of a special case of the Huang-Ji theorem that real ellipsoids possess umbilical points; namely, we prove the existence of real curves of umbilical points in the special case where the ellipsoids are sufficiently close to spherical. This is Theorem \ref{Prop-ell}.

One of the main obstacles in investigating the existence of umbilical points on a compact real hypersurface $M^3\subset\bC^2$ is the lack of a global convenient representation ("formula") of Cartan's tensor $Q$. In Chern--Moser normal form at $p=(0,0)$, $Q$ is represented at $p$ by the coefficient $c_{24}$. Loboda \cite{Loboda97} discovered a "semi-global" formula that represents $Q$ for a graph $M^3\subset \bC^2$ provided that $M^3$ also has additional tranverse symmetry. One of the main results in this paper is a {\it global} representation of $Q$ as a nonlinear PDO acting on a defining function $\rho$ for $M^3\subset \bC^2$. In fact, a sequence of invariants $\det A_k(\rho)$, for $k\geq 3$, is introduced in Section \ref{SecInv} (see in particular Theorem \ref{transform}). These invariant determinants are in some sense higher order versions of Fefferman's Monge-Ampere operator $J(\rho)$. If we introduce the $(1,0)$-vector field
$$
L:=-\frac{\partial \rho}{\partial w}\frac{\partial}{\partial z}+\frac{\partial \rho}{\partial z}\frac{\partial}{\partial w}=-\rho_w\partial_z+\rho_z\partial_w,
$$
then we have, as the reader can easily check,
\beq\Label{FeffJ-2}
J(\rho)=\det
\begin{pmatrix}
\rho_z & \bar L\rho_z\\
\rho_w & \bar L \rho_w
\end{pmatrix} \mod \rho.
\eeq
Our invariant determinants $\det A_k(\rho)$, $k\geq 3$ are higher order analogues of the determinant on the right in \eqref{FeffJ-2}. The precise definition is given in \eqref{an}. The relationship between $\det A_k(\rho)$ and the classical Chern--Moser invariants of $M^3$ is explored in Section \ref{SecCM}. This relationship is expressed by \eqref{detAnCM}; in particular, it is shown that $\det A_3(\rho)$ represents Cartan's tensor $Q$ at every $p\in M^3$.

\section{A family of invariant determinants}\Label{SecInv}

Consider a real hypersurface $M=M^3\subset\C^2$ given by $\rho(z,w,\bar z,\bar w)=0$
with $\d\rho\ne 0$ on $M$. We use coordinates $Z=(z,w)\in \C^2$.
Then the space of $(1,0)$-vectors on $M$ at every point is spanned by
\beq\Label{l}
L: = - \rho_w \d_z + \rho_z \d_w.
\eeq
We shall also use the Hessian $\rho_{Z^2}$ evaluated at $(L,L)$:
$$\rho_{Z^2}(L,L) = \rho_{zz} \rho_2^2 - 2\rho_{zw} \rho_z\rho_w + \rho_{w^2} \rho_z^2. $$
For every $n\ge 3$, consider the $(2n-1)\times(2n-1)$ matrix PDO $A_n(\rho)$ acting on the smooth function $\rho$
\beq\Label{an}
A_n=A_n(\rho):=
\begin{pmatrix}
\bar L^j (\rho_z^k \rho_w^{n-k})\cr
\bar L^j (\rho_z^s \rho_w^{n-3-s} \rho_{Z^2}(L, L))
\end{pmatrix}_{0\le j\le 2n-2,\, 0\le k\le n, \, 0\le s\le n-3},
\eeq
where we regard $j$ as column index and $k$ and $s$ as row indices (first followed by the second). In particular, for $n=3,4$, we obtain
\beq\Label{a3}
A_3=A_3(\rho):=
\begin{pmatrix}
\rho_w^3  & \bar L(\rho_w^3) & \cdots & \bar L^4(\rho_w^3)\cr
 \rho_z\rho_w^2 & \bar L( \rho_z\rho_w^2) & \cdots & \bar L^4( \rho_z\rho_w^2)\cr
  \rho_z^2\rho_w & \bar L(  \rho_z^2\rho_w) & \cdots & \bar L^4(  \rho_z^2\rho_w)\cr
\rho_z^3 & \bar L(\rho_z^3) & \cdots & \bar L^4(\rho_z^3)\cr
\rho_{Z^2}(L, L) & \bar L(\rho_{Z^2}(L, L)) & \cdots & \bar L^4(\rho_{Z^2}(L, L))
\end{pmatrix}
\eeq
and
\beq\Label{a4}
A_4=A_4(\rho):=
\begin{pmatrix}
\rho_w^4  & \bar L(\rho_w^4) & \cdots & \bar L^6(\rho_w^4)\cr
 \rho_z\rho_w^3 & \bar L( \rho_z\rho_w^3) & \cdots & \bar L^6( \rho_z\rho_w^3)\cr
\rho_z^2\rho_w^2  & \bar L(\rho_z^2\rho_w^2) & \cdots & \bar L^6(\rho_z^2\rho_w^2)\cr
 \rho_z^3\rho_w  & \bar L( \rho_z^3\rho_w ) & \cdots & \bar L^6( \rho_z^3\rho_w )\cr
\rho_z^4 & \bar L(\rho_z^4) & \cdots & \bar L^6(\rho_z^4)\cr
\rho_w\rho_{Z^2}(L, L) & \bar L(\rho_w\rho_{Z^2}(L, L)) & \cdots & \bar L^6(\rho_w\rho_{Z^2}(L, L))\cr
\rho_z\rho_{Z^2}(L, L) & \bar L(\rho_z\rho_{Z^2}(L, L)) & \cdots & \bar L^6(\rho_z\rho_{Z^2}(L, L))
\end{pmatrix}.
\eeq
We also denote by $D_n=D_n(\rho)$ the upper left $(n+1)\times (n+1)$ minor of $A_n$, i.e.\
\beq\Label{dn}
D_n:=
\begin{pmatrix}
\bar L^j (\rho_z^k \rho_w^{n-k})
\end{pmatrix}_{0\le j\le n, \, 0\le k\le n}.
\eeq
The main interest in considering the matrices $A_n$ and $D_n$ is the following invariance property of their determinants:

\bt\Label{transform}
For every $M$ and $n\ge 3$, the properties $\det A_n=0$ and $\det D_n=0$ at points of $M$
are independent of the choice of the defining function $\rho$ as well as of
the choice of the coordinates $Z=(z,w)\in \C^2$.

More precisely, if $L^*$, $A^*_n$ and $D^*_n$ are given by \eqref{l}, \eqref{an} and \eqref{dn} respectively
with $\rho$ replaced by another defining function $\rho^*=a\rho$ (where $a$ is any nonzero real smooth function),
and $Z=(z,w)$ replaced by another (formal) holomorphic coordinate system $Z^*=(z^*,w^*)$,
we have the transformation rule
\beq
\Label{a-transform}
\delta^{n^2-1} \bar \delta^{(n-1)(2n-1)}\det A^*_n = a^{(2n-1)^2}  \det A_n,
\eeq
\beq
\Label{d-transform}
|\delta|^{n(n+1)} \det D^*_n = a^{\frac{3n(n+1)}2}  \det D_n,
\eeq
where $\delta$ is the Jacobian determinant of the coordinate transformation  $Z^*=H(Z)$.
\et

\br
It is important that $L$ and $\rho$ used in $A_n$ are related via \eqref{l}.
The invariance of the property $\det A_n=0$
does not hold for arbitrary choices of $(1,0)$ vector fields $L$.
However, given that chosen $L$, the invariance of $\det A_n$ remains
when replacing $\bar L$ by arbitrary $(0,1)$-vector fields.
\er

To prove Theorem \ref{transform}, we require two lemmas.

\bl\Label{a-lemma}
For any real smooth function $a=a(Z,\bar Z)$, and $L$ and $L^*$ given respectively by \eqref{l} and by the same formula with $\rho$ replaced with $\rho^*=a\rho$, we have on $M$ the identities
\beq\Label{a-transform'}
\rho^*_Z = a\rho_Z,
\quad L^* = a L,
\quad \rho^*_{Z^2}(L^*,L^*) = a^3 \rho_{Z^2}(L, L).
\eeq
\el

\bpf
By Leibnitz' rule on $M$ we have $\rho^*_Z = a \rho_Z $, which implies the first and second identities in \eqref{a-transform'}. By Leibnitz' rule again, for any $(1,0)$ vectors $\xi, \eta$, we also have
$$(a\rho)_{Z^2}(\xi, \eta)
= a\rho_{Z^2}
(\xi,\eta)
+ a_Z(\xi) \rho_Z(\eta)
+ a_Z(\eta) \rho_Z(\xi)
+ a_{Z^2}(\xi,\eta) \rho.
$$
Substituting $\xi=\eta=L$ and using the properties $\rho_Z(L)=0$ and $\rho=0$ on $M$,
we obtain, on $M$,
\beq\Label{rho-L}
\rho^*_{Z^2}(L, L) = a \rho_{Z^2} (L, L).
\eeq
Together with the second identity in \eqref{a-transform'} this yields the third identity.
\epf

\bl\Label{coord}
For any (formal) biholomorphic transformations $Z^*=H(Z)$, $\zeta^*=K(\zeta)$ of $\C^2$, any complex formal power series $\rho^*(Z^*, \bar Z^*)$
and $\rho(Z,\bar Z):= \rho^*(H(Z), K(\bar Z))$,
consider $L$, $A_n$, $D_n$ and $L^*$, $A^*_n$, $D^*_n$ given by \eqref{l}, \eqref{an}, \eqref{dn} and respectively
by the same identities with $\rho$ replaced by $\rho^*$.
Then the following hold:
\begin{enumerate}
\item[(i)] The identity
\beq\Label{Z2}
\rho_{Z^2}(L,L) =  (\det H_Z)^2 \rho^*_{Z^{*2}} (L^*, L^*)
\eeq
holds modulo a cubic homogeneous polynomial in $(\rho_z, \rho_w)$ with holomorphic coefficients in $Z$.
\item[(ii)] The determinants of $A_n$, $D_n$ and $A^*_n$, $D^*_n$ are related by
\beq
\det A_n = (\det H_Z)^{n^2-1}(\det K_{\bar Z})^{(n-1)(2n-1)} \det A^*_n,
\eeq
\beq
\det D_n = (\det H_Z)^{\frac{n(n+1)}2}
(\det K_{\bar Z})^{\frac{n(n+1)}2} \det D^*_n,
\eeq
where the matrices $A^*_n$ and $D^*_n$ are evaluated at $(Z^*,\bar Z^*)= (H(Z), K(\bar Z))$.
\end{enumerate}
\el

\bpf
We write
$$S_n:=
\begin{pmatrix}
\rho_w^n\cr
\rho_z\rho_w^{n-1}\cr
\vdots\cr
\rho_z^n
\end{pmatrix}
$$
for the standard basis of homogeneous monomials of order $n$ in $(\rho_z, \rho_w)$.
In particular, $S_n$ coincides with the
 $(n+1)\times 1$ matrix (column vector) consisting of the first $n+1$ entries of the first column of $A_n$.
We also write $S^*_n$ for corresponding column of monomials in $(\rho^*_{z^*}, \rho^*_{w^*})$.
In particular,
$$
S_1=
\begin{pmatrix}
\rho_w\cr
\rho_z
\end{pmatrix},
\quad
S^*_1=
\begin{pmatrix}
\rho^*_{w^*}\cr
\rho^*_{z^*}
\end{pmatrix}.
$$
By the chain rule, we have
$$\rho_Z = \rho^*_{Z^*}\circ H_Z,$$
where $\rho^*_{Z^*}$ is evaluated at $(Z^*,\bar Z^*)=(H(Z), K(\bar Z))$.
Writing as matrix identity we obtain
$$S_1 = H_{Z} S^*_1,$$
where by abuse of notation, for $H=(f,g)$, we identify $H_Z$ with its induced matrix
\beq\Label{h-matrix}
\begin{pmatrix}
g_w & f_w\cr
g_z & f_w
\end{pmatrix}.
\eeq
Then viewing $n$th order homogenous monomials in $\rho_z,\rho_w$ in the $n$th tensor power of the cotangent space, we have
\beq\Label{S-rel}
S_n = (\otimes^n H_Z)  S^*_n,
\eeq
where, by another abuse of notation, we identify the tensor power transformation $\otimes^n H_Z$
with its induced matrix in the monomial basis.
Further we have
\beq\Label{det}
\det (\otimes^n H_Z) = (\det H_Z)^{\frac{n(n+1)}2},
\eeq
which e.g. follows from Jordan normal form.
Since $(\otimes^n H_Z)$ is holomorphic in $Z$ and $\bar L$ is a $(0,1)$ vector field, we obtain for any $j$,
\beq\Label{first-rows}
\bar L^j S_n =  (\otimes^n H_Z)  \bar L^j S^*_n.
\eeq
In particular, we have
\beq\Label{dn-rel}
\det D_n = (\det H_Z)^{\frac{n(n+1)}2} \det E^*_n,
\eeq
where $E^*_n$ is the matrix obtained from $D^*_n$ with $\bar L^*$
being replaced by $\bar L$.

We next turn to the relation between $L$ and $L^*$.
By definition
$$
L=
\begin{pmatrix}
\rho_z & -\rho_w
\end{pmatrix}
\begin{pmatrix}
\d_w\cr
\d_z
\end{pmatrix},
\quad
L^*=
\begin{pmatrix}
\rho^*_{z^*} & -\rho^*_{w^*}
\end{pmatrix}
\begin{pmatrix}
\d_{w^*}\cr
\d_{z^*}
\end{pmatrix}.
$$
We further have
$$
\begin{pmatrix}
H_Z \d_w \cr
H_Z \d_z
\end{pmatrix}
=
\begin{pmatrix}
g_w & f_w\cr
g_z & f_z
\end{pmatrix}
\begin{pmatrix}
\d_{w^*} \cr
\d_{z^*}
\end{pmatrix}
= H_Z
\begin{pmatrix}
\d_{w^*} \cr
\d_{z^*}
\end{pmatrix},
$$
where we continue our abuse of notation by writing $H_Z$ also for its induced matrix.
Similarly,
\beq\Label{rho-rel}
\begin{pmatrix}
\rho_w & \rho_z
\end{pmatrix}
=
\begin{pmatrix}
\rho^*_{w^*} & \rho^*_{z^*}
\end{pmatrix}
\begin{pmatrix}
g_w & g_z \cr
f_w & f_z
\end{pmatrix}
=
\begin{pmatrix}
\rho^*_{w^*} & \rho^*_{z^*}
\end{pmatrix}
H_Z^t,
\eeq
where $H_Z^t$ is the transpose matrix.
Furthermore,
$$
\begin{pmatrix}
\rho_z & -\rho_w
\end{pmatrix}
=
\begin{pmatrix}
\rho_w & \rho_z
\end{pmatrix}
J,
\quad
J:=
\begin{pmatrix}
0 & -1\cr
1 &0
\end{pmatrix},
$$
and hence
$$
\begin{pmatrix}
\rho_z & -\rho_w
\end{pmatrix}
=
\begin{pmatrix}
\rho^*_{z^*} & -\rho^*_{w^*}
\end{pmatrix}
J^{-1} H_Z^t J.
$$
Putting everything together, we obtain
$$H_Z L =
\begin{pmatrix}
\rho_z & -\rho_w
\end{pmatrix}
\begin{pmatrix}
H_Z \d_w \cr
H_Z \d_z
\end{pmatrix}
=
\begin{pmatrix}
\rho_z & -\rho_w
\end{pmatrix}
H_Z
\begin{pmatrix}
\d_{w^*} \cr
\d_{z^*}
\end{pmatrix}
=
\begin{pmatrix}
\rho^*_{z^*} & -\rho^*_{w^*}
\end{pmatrix}
J^{-1} H_Z^t J H_Z
\begin{pmatrix}
\d_{w^*} \cr
\d_{z^*}
\end{pmatrix}.
$$
By direct calculation,
$$
J^{-1} H_Z^t J H_Z =
\begin{pmatrix}
0 & 1\cr
-1 &0
\end{pmatrix}
\begin{pmatrix}
g_w & g_z \cr
f_w & f_z
\end{pmatrix}
\begin{pmatrix}
0 & -1\cr
1 &0
\end{pmatrix}
\begin{pmatrix}
g_w & f_w \cr
g_z & f_z
\end{pmatrix}
= (\det H_Z)
\begin{pmatrix}
1 & 0\cr
0 &1
\end{pmatrix},
$$
and therefore
\beq\Label{LL}
H_Z L = (\det H_Z) L^*.
\eeq
Now, by the chain rule, for any $(1,0)$ vectors $\xi, \eta$,
$$
\rho_{Z^2}(\xi,\eta) = \rho^*_{Z^{*2}} (H_Z \xi, H_Z \eta) + \rho^*_{Z^*} (H_{Z^2}(\xi, \eta)),
$$
where we recall that $H_{Z^2}(\xi, \eta)$ is a $(1,0)$ vector.
Substituting $\xi=\eta = L$ and using \eqref{LL}, we obtain
$$
\rho_{Z^2}(L,L) =  (\det H_Z)^2 \rho^*_{Z^{*2}} (L^*, L^*) + \rho^*_{Z^*}( H_{Z^2}(L,L)).
$$
Expanding the last term and using \eqref{rho-rel}, we obtain the desired identity \eqref{Z2}
modulo a homogeneous polynomial of degree 3 in $(\rho_z,\rho_w)$ with coefficients that are polynomial in the derivatives of $H$, which proves the first statement (i) of the lemma.

We now turn to the $(n-2)\times 1$ matrix (column vector) consisting of the last $n-2$ entries of the first column in $A_n$, which in the notation introduced can be expressed as $\rho_{Z^2}(L,L)S_{n-3}$.
Hence \eqref{S-rel} implies
$$\rho_{Z^2}(L,L) S_{n-3} = \rho_{Z^2}(L,L) (\otimes^{n-3} H_Z)  S^*_{n-3}.$$
By the first statement (i) of the lemma, already proved, we have
$$\rho_{Z^2}(L,L) S_{n-3} = (\det H_Z)^2  (\otimes^{n-3} H_Z)  \rho^*_{Z^{*2}}(L^*,L^*) S^*_{n-3}$$
modulo a homogeneous polynomial of order $n$ in $(\rho_z, \rho_w)$ with holomorphic coefficients in $Z$.
Since $\bar L$ is $(0,1)$, it commutes with those holomorphic coefficients and, consequently,
we can subtract from the last $n-2$ rows of $A_n$ suitable linear combinations of the first $n+1$ rows
to obtain a matrix with the same determinant of the form
$$
\begin{pmatrix}
\otimes^n H_Z & 0\cr
0 & (\det H_Z)^2 (\otimes^{n-3} H_Z)
\end{pmatrix}
B^*_n,
$$
where $B^*_n$ is the matrix obtained from $A^*_n$ with $\bar L^*$ (but not $L^*$!) replaced by $\bar L$.

Finally, analogously to \eqref{LL}, we have the relation
$$
K_{\bar Z}\bar L = (\det K_{\bar Z}) \bar L^*.
$$
Writing $C^*_n$ for the first column of $B^*_n$,
we obtain for any $j$,
$$(K_{\bar Z}\bar L)^j C^*_n =   (\det K_{\bar Z})^j (\bar L^*)^j C^*_n$$
modulo a linear combination of the columns $(\bar L^*)^s C^*_n$ with $s<j$.
Hence, subtracting those linear combinations without changing the determinant, we obtain
$$\det B^*_n = (\det K_{\bar Z})^{(n-1)(2n-1)} \det A^*_n, $$
and similar
$$\det E^*_n = (\det K_{\bar Z})^{\frac{n(n+1)}2} \det D^*_n,$$
where $E^*_n$ was defined after \eqref{dn-rel}.
Putting everything together and using \eqref{det} we obtain the second conclusion.
\epf

\bpf[Proof of Theorem~$\ref{transform}$]
Clearly it suffices to prove the proposition by separately considering changes of the defining function and the coordinates. The transformation formula under a change of coordinates follows from Lemma~\ref{coord} when $\rho$ is real-analytic. However, since the matrix $A_n$
only depends on a finite order jet of $\rho$ at a reference point,
the corresponding transformation rule in \eqref{a-transform'} holds for any smooth $\rho$.

It remains to consider the change $\rho^* = a\rho$ of the defining function.
By Lemma~\ref{a-lemma}, for any $k$, $s$, $n$ as in \eqref{an}, we have
$$(\rho^*_z)^k(\rho^*_w)^{n-k} = a^n \rho_z^k \rho_w^{n-k},
\quad
(\rho^*_z)^s(\rho^*_w)^{n-3-k} \rho^*_{Z^2}(L^*, L^*)
 = a^n \rho_z^k \rho_w^{n-k} \rho_{Z^2}(L, L).
$$
Then, writing
$$
C_n:=
\begin{pmatrix}
(\rho_z^k \rho_w^{n-k})\cr
(\rho_z^s \rho_w^{n-3-s} \rho_{Z^2}(L, L)
\end{pmatrix}_{0\le k\le n, \, 0\le s\le n-3}
$$
for the first column of the matrix $A_n$ given by \eqref{an}
and $C^*_n$ for the corresponding first column of $A^*_n$, we obtain
$$C^*_n = a^n C_n.$$
Then using the relation $\bar L^* = a\bar L$, we conclude for every $j$,
$$(\bar L^*)^j C^*_n = a^{n+j} \bar L^j C_n$$
modulo linear combinations of the columns $\bar L^s C_n$ with $s<j$.
Since the determinant does not change after subtracting a linear combination for columns from another column, we obtain the desired transformation rules
$$\det A^*_n = a^{(2n-1)^2} \det A_n,
\quad
\det D^*_n = a^{\frac{3n(n+1)}2} \det D_n.
$$
\epf

\section{Calculation in Chern-Moser normal form}\Label{SecCM}

Note that the invariance property in Proposition~\ref{det} was obtained
for any smooth real hypersurface given by $\rho(Z,\bar Z)=0$.
If the latter is Levi-nondegenerate, we can use special defining functions
\beq\Label{normform}
\rho(z,w,\bar z,\bar w) = -\Im w + \phi(z,\bar z, \Re w), \quad \phi(z,\bar z, u) = \sum \phi_{kl}(u)z^k \bar z^l,
\eeq
in Chern-Moser normal form (in the formal sense if $M$ is only smooth and not real-analytic) to
compute the determinant of $A_n$ at the origin.
Recall \cite{CM74} that the normal form requires $\phi$ to satisfy
\beq\Label{CMnormal}
\phi_{11}=1, \quad
\phi_{0k}=\phi_{1s}=\phi_{22}=\phi_{23}=\phi_{33}=0, \quad k\ge 0, s\ge 2.
\eeq
In this normal form we have $(\rho_w, \rho_z)(0)=(i/2,0)$
and furthermore
$$\bar L^j \rho_w(0)= \bar L^k \rho_z(0) = 0, \quad j\geq 1,\ k\ne 1,$$
and
$$\bar L\rho_z(0) =-i/2 \ne 0. $$
Furthermore,
$$\rho_{Z^2}(L,L) = \rho_w^2 \rho_{z^2} - 2 \rho_z\rho_w \rho_{zw} + \rho_{z}^2 \rho_{w^2}  $$
and
$$\rho_{\bar z^l w^s}(0) = \rho_{z \bar z^l w^s}(0) = 0, \quad l\ge 0, \quad s\ge 1,$$
imply
$$\bar L^k (\rho_{Z^2}(L,L))(0) =
(\rho_w(0))^2 \bar L^k \rho_{z^2}(0) =
 (-i/2)^{k+2} \rho_{z^2\bar z^k} (0) =
 (-i/2)^{k+2} \phi_{z^2\bar z^k} (0).$$
We similarly observe that
\begin{multline}
\bar L^k (\rho_z^s\rho_w^{n-3-s}\rho_{Z^2}(L,L))(0) =\binom{k}{s}
(\rho_w(0))^{n-1-s}\bar L^s\rho_z^s (0)\bar L^{k-s} \rho_{z^2}(0) \\= (-1)^k
 (i/2)^{n+k-s-1} s!\binom{k}{s}\phi_{z^2\bar z^{k-s}} (0).
 \end{multline}
Hence we obtain
\beq\Label{detAnCM}
\det A_n|_{Z=0} = c_n \det
\begin{pmatrix}
\binom{n+1}{0}\phi_{2,n+1} & \binom{n+2}{0}\phi_{2,n+2} & \cdots & \binom{2n-2}{0}\phi_{2, 2n-2}\cr
\binom{n+1}{1}\phi_{2,n} & \binom{n+2}{1}\phi_{2,n+1} & \cdots & \binom{2n-2}{1}\phi_{2, 2n-3}\cr
\vdots & \vdots & \ddots &\vdots\cr
\binom{n+1}{n-3}\phi_{2,4} & \binom{n+2}{n-3}\phi_{2, 5} & \cdots & \binom{2n-2}{n-3}\phi_{2, n+1}\cr
\end{pmatrix},
\eeq
where $c_n\ne 0$ is a universal constant (independent of $\phi$).
In particular,
\beq\Label{A3=Q}
\det A_3|_{Z=0} = c_3 \phi_{2,4},\eeq
is Cartan's ``6th order tensor",
$$\det A_4|_{Z=0} = c_4 \det
\begin{pmatrix}
\phi_{2,5} & \phi_{2,6}\cr
5\phi_{2, 4} & 6\phi_{2,5}
\end{pmatrix},
$$
$$\det A_5|_{Z=0} = c_5 \det
\begin{pmatrix}
\phi_{2,6} & \phi_{2,7} & \phi_{2,8}\cr
6\phi_{2, 5} & 7\phi_{2,6} & 8\phi_{2,7}\cr
15\phi_{2, 4} & 21\phi_{2,5} & 28\phi_{2,6}
\end{pmatrix}.
$$

The same calculations, for any hypersurface in the form \eqref{normform} with $\phi(z,0,s)=\phi(0,\bar z,s)\equiv 0$ as a formal power series (which can always be achieved; see \cite{BER99a}), show also
that each $\det D_n$ (given by \eqref{dn}) equals a universal constant times $(\phi_{11})^{n}$.
Together with Theorem~\ref{transform} this yields:

\bp\Label{dn-nonvanish} Let $M\subset \bC^2$ be a real smooth hypersurface
$M$ given by $\rho(Z,\bar Z)=0$. Then, $\det D_n\ne 0$ at $p\in M$ if and only if
$M$ is Levi-nondegenerate at $p$.
\ep

\section{Umbilical points on real hypersurfaces in $\bC^2$}

Let $M\subset \bC^2$ be a smooth real hypersurface and $p\in M$. Recall that $p$ is said to be an {\it umbilical point} if in Chern--Moser normal coordinates $(z,w)$ vanishing at $p$, the coefficient $\phi_{2,4}$ in the Chern--Moser normal form (\eqref{normform} and \eqref{CMnormal}) vanishes, $\phi_{2,4}=0$. While Chern-Moser normal coordinates and normal form are not unique, it is well known \cite{CM74} that the vanishing of $\phi_{2,4}$ is an invariant. By Theorem \ref{transform} and \eqref{A3=Q}, we have:

\bp\Label{Umb-A3} Let $M\subset \bC^2$ be defined by $\rho=0$. Then, the set $\mathcal U$ of umbilical points on $M$ is given by the equation $\det A_3(\rho)=0$.
\ep

\subsection{Umbilical indices}
For a fixed global defining equation $\rho=0$ for $M$, where $\rho$ is defined in a neighborhood of $M$, denote by $Q:=\det A_3(\rho)$, so that the set $\mathcal U\subset M$ of umbilical points is given by $Q=0$. We assume $M$ to be oriented and choose $\rho$ compatible with that orientation, i.e.\ such that the gradient of $\rho$ completes positively oriented frames in $M$ to those in $\C^2$. Note that such $\rho$ always exists e.g.\ the oriented distance function. Further, $\rho$ is unique up to multiplication with a positive real function in a neighborhood of $M$ in $\C^2$.

We shall say that $p\in \mathcal U\subset M$ is a {\it $1$-regular} umbilical point of $M$ if $\mathcal U$ is a smooth real curve (1-manifold) at $p$. By Thom's transversality, every hypersurface $M$
can be approximated by one having only $1$-regular umbilical points.

\bd
For every oriented closed curve $C$ in $M$ avoiding the umbilic set $\mathcal U$,
define its {\em umbilical index} to be $-1/2$ times the winding number
of $Q$ along $C$. For every $1$-regular umbilical point $p$
with chosen orientation of $\mathcal U$,
define its {\em local umbilical index}
(or umbilical index of $M$ at $p$)
to be the umbilical index of the positively oriented boundary
of any sufficiently small disk transversal to $\mathcal U$.
\ed

Since $\rho$ is unique up to multiplication with positive real function,
it follows from Theorem~\ref{transform} that the index as defined
is independent of the choice of $\rho$. It further follows
from the same theorem that the umbilic index of $C$
is also independent of the choice of the ambient coordinates in $\C^2$
as long as $C$ is null-homotopic.

If $\Sigma\subset M$ is an oriented surface (2-manifold) that meets $\mathcal U$ transversely at a $1$-regular point $p$,
the index of $M$ at $p$ is given by
\beq\Label{Indexp}
\iota_\Sigma(p)=-\frac{1}{2}\deg\left(\frac{Q}{|Q|}\colon \partial \Sigma_p\to S^1\right),
\eeq
where $\Sigma_p$ is the boundary of
 a sufficiently small topological disk containing $p$
  (topologically a circle $S^1$) oriented positively with respect to $\Sigma$, and where $\deg$ denotes the topological degree
  (which is the same as winding number in this case).
 We note that we can also express the index using an integral,
\beq\Label{Indexp}
\iota_\Sigma(p)=\frac{i}{4\pi }\int_{\partial\Sigma_p}\frac{dQ}{Q}=
\frac{i}{4\pi }\int_{\partial\Sigma_p}d\log Q.
\eeq
Since the 1-form $dQ/Q$ is closed away from the zeros of $Q$, we observe from \eqref{Indexp} that in fact the index at $p$ only depends on the orientation of $\Sigma$ and not on the choice of transversal $\Sigma$ itself, and that the index is constant along every component of the 1-manifold of nondegenerate umbilical points $\mathcal U_{1}$. If $p\in M$ is a $1$-regular umbilical point and the index of $M$ at $p$ (with respect to any transversal $\Sigma$) is non-zero, then we shall say that $p$ is a {\it stable} umbilical point. In view of Thom's transversality, any sufficiently small perturbation of the CR structure of $M$ near a stable umbilical point $p$ will have stable umbilical points near $p$, which motivates this terminology.

We shall use the notation $W_\gamma(R)$ for the winding number of a function $R$ on $M$ along an oriented closed curve $\gamma$ (defined only when $R$ does not vanish on $\gamma$);
\beq\Label{W(R)}
W_\gamma(R)=\frac{1}{2\pi i}\int_{\gamma}\frac{dR}{R}=
\frac{1}{2\pi i}\int_{\gamma}d\log R.
\eeq
Thus, in particular, if $\Sigma$ is a surface, transversal to $\mU$ at $p$ and $\Sigma_p$ is its intersection with a small tubular neighborhood of $\mu$ near $p$, then by definition of the index:
$$
\iota_\Sigma(p)=-\frac{1}{2}W_{\partial\Sigma_p}(Q).
$$

If $M$ is real-analytic, then $\mathcal U$ is a real-analytic subvariety of $M$, and the set $\mathcal U_1$ of 1-regular umbilical points consist of the subset of $\mathcal U$ of regular points of dimension one. For simplicity, we shall proceed under the assumption that $M$ is real-analytic. In this case, $\mathcal U$ is either all of $M$ (we assume that $M$ is connected), in which case $M$ is locally spherical, or $\mathcal U$ is a proper subvariety. In the latter case, points of $\mathcal U$ are either 0-, 1-, or, 2-dimensional, and we decompose $\mathcal U$ accordingly, $\mathcal U=\mathcal U^0\cup \mathcal U^1\cup \mathcal U^2$; recall that the (topological) dimension of a real-analytic subvariety $\mathcal V$  at a point $p$ is the largest dimension of a nonsingular component of $\mathcal V$ with $p$ in its closure.

We note that the set of $1$-regular umbilical points $\mU_1$ equals $\mathcal U^0\cup \mathcal U^1$ minus a discrete set of points. Thus, if there are no points of dimension $2$ on $\mU$, then any surface $\Sigma$ that intersects $\mU$ can be locally perturbed to only intersect $\mU$ along the set of $1$-regular points $\mU_1$. We have the following simple consequence of Stokes Theorem:

\bp\Label{Wvsindex}
Let $M\subset \bC^2$ be a real-analytic hypersurface, and assume that $\mathcal U$ has no points of dimension $2$ or $3$, i.e., $\mathcal U=\mU^0\cup\mU^1$. Let $\gamma\subset M$ be a oriented closed curve, homologous to $0$ in $M$ and not intersecting $\mU$, and $\Sigma\subset M$ an oriented surface, intersecting $\mU$ transversally along $\mU_1$ and with $\partial\Sigma=\gamma$. Then,
\beq
W_\gamma(Q)=-\frac{1}{2}\sum_{p\in \mU_1\cap \Sigma}\iota_\Sigma(p).
\eeq
\ep

\begin{proof} Let $p_1,\ldots, p_k$ be the finite set of points in $\mU_1\cap \Sigma$ and $\Sigma_{p_j}$ the intersection of $\Sigma$ with a sufficiently small tubular neighborhood of $\mU_1$ near $p_j$ (so small that the closures of the $\Sigma_{p_j}$ are disjoint). Since $dQ/Q$ is closed in the punctured surface
$$
\Sigma':=\Sigma\setminus \bigcup_{j=1}^k \Sigma_{p_j},
$$
being locally $d\log Q$ on $\Sigma'$, we conclude by Stokes Theorem that:
\beq
W_\gamma(Q)=\frac{1}{2\pi i} \int_{\partial\Sigma}\frac{dQ}{Q}=\sum_{j=1}^k \frac{1}{2\pi i}\int_{\partial\Sigma_{p_j}} \frac{dQ}{Q}= -\frac{1}{2}\sum_{p\in \mU_1\cap \Sigma}\iota_\Sigma(p).
\eeq
\end{proof}

\section{Umbilical points on Real Ellipsoids}\Label{EllSec}

We shall consider real ellipsoids $E\subset \bC^2$. A general real ellipsoid can be defined by an equation of the form
\begin{equation}
A(z^2+\bar z^2)+2(2+A)|z|^2+B(w^2+\bar w^2)+2(2+B)|w|^2=4,\quad A,B\geq 0.
\end{equation}
We shall fix $A,B\geq 0$, and assume that at least one of these, say $A$, is nonzero (so that the ellipsoid does not degenerate to a sphere). We consider a 1-parameter family of ellipsoids $E_\epsilon$, defined by $\rho_\epsilon=0$, where
\begin{equation}
\rho_\epsilon:=
\epsilon A(z^2+\bar z^2)+2(2+\epsilon A)|z|^2
+\epsilon B(w^2+\bar w^2)+2(2+\epsilon B)|w|^2
-4
\end{equation}
\beq
= -4 + 4(|z|^2 + |w|^2) + \epsilon (
	A (z^2 + \bar z^2 + 2|z|^2 )
	+ B (z^2 + \bar z^2 + 2|z|^2
	),
	\quad \epsilon >0.
\eeq
Note that $E_0$ is the unit sphere. We shall mainly be concerned with small perturbations of the sphere, and shall thus consider small $\epsilon>0$. Our aim is to prove the following result:

\begin{Thm}\Label{Prop-ell} For $\epsilon>0$ sufficiently small, the subset of umbilical points on $E_\epsilon$ either contains points of dimension at least $2$ or contains a curve of umbilical points.
\end{Thm}

To this end, we shall compute the matrix $A_3=A_3(\rho_\epsilon)$ on $\rho=\rho_\epsilon=0$. Since we will be concerned with small perturbations it suffices, as we shall see, to compute $A_3$ mod $O(\epsilon^3)$. We note first that
\begin{equation}\label{L-ell}
\rho_z=4\bar z+2\epsilon A (\bar z +z),\quad \rho_w=4\bar w+2\epsilon B(\bar w+ w)
\end{equation}
and
\beq
	\rho_{z^2} = 2\epsilon A
	,\quad \rho_{zw} = 0
	,\quad \rho_{w^2} = 2\epsilon B.
\eeq
Therefore, we have
\begin{equation}\label{rhoZ2}
\begin{aligned}
\rho_{Z^2}(L,L)=&(\rho_w)^2\rho_{z^2}-2\rho_z\rho_w\rho_{zw}+(\rho_z)^2\rho_{w^2}\\
=&8\epsilon \left(A(2\bar w+\epsilon B(\bar w+ w))^2+B(2\bar z+\epsilon A (\bar z +z))^2\right)\\
=&32\epsilon (A\bar w^2+B\bar z^2)
	 + 32\epsilon^2AB(|z|^2+|w|^2+\bar w^2+\bar z^2)\\
	& + 8\epsilon^3AB\left (A(z+\bar z)^2+B(w+\bar w)^2\right).
\end{aligned}
\end{equation}
To calculate $A_3(\rho)$, we shall repeatedly apply $\bar L$, where
\begin{equation}\label{barL}
\begin{aligned}
	\bar L &=
		-\rho_{\bar w} \frac{\partial}{\partial \bar z}
		+ \rho_{\bar z}\frac{\partial}{\partial \bar w}\ \\
	&=-2(2w+\epsilon B(\bar w+ w))\frac{\partial}{\partial \bar z}+2(2z+\epsilon A(\bar z+z)\frac{\partial}{\partial \bar w}\\
&=4\left(-w\frac{\partial}{\partial \bar z}+z\frac{\partial}{\partial \bar w}\right )+2\epsilon\left(-B(w+\bar w)\frac{\partial}{\partial \bar z}+A(z+\bar z)\frac{\partial}{\partial \bar w}\right)\\
&=\bar L_0+\epsilon \bar L_1,
\end{aligned}
\end{equation}
to $\rho_{Z^2}(L,L)$,
and subsequently evaluate at $w=0$. Since $\bar L$ only involves differentiation in $\bar z$ and $\bar w$,
the result for $w=0$ will not change when replacing with $0$
 all occurences of
$w$ (but not $\bar w$). Thus for $w=0$, we obtain
\beq
\begin{aligned}
& \frac1{2^k \cdot 32}  \bar L^k \rho_{Z^2}(L,L)  \\
& =  \left(
	2z \frac\d{\d\bar w}
	+ \epsilon \left(
		- B\bar w \frac\d{\d\bar z}
		+ A(z + \bar z) \frac\d{\d\bar w}
	\right)
\right)^k
\left(
	\epsilon (A \bar w^2 + B \bar z^2)
	+ \epsilon^2 AB(|z|^2 + \bar z^2 + \bar w^2).
\right)
\end{aligned}
\eeq
Then we obtain for $w=0$,
\beq
	\bar L\rho_{Z^2}(L,L)
	= \bar L^3 \rho_{Z^2}(L,L)
	= \bar L^4 \rho_{Z^2}(L,L)
	= O(\epsilon^2),
	\quad
	\frac1{2^2 \cdot 32} \bar L^2\rho_{Z^2}(L,L)
	= 2^3 z^2 A \epsilon + O(\epsilon^2)
\eeq

\subsection{Terms of the form $\bar L^k\rho_z^3$; first row.} We note that
\beq
\rho_z^3=(4\bar z+2\epsilon A (\bar z +z))^3=8(2\bar z+\eps A(\bar z+z))^3.
\eeq
We shall be interested in $A_3$ mod $O(\epsilon^3)$, and since all terms in the last row (computed above) are already $O(\eps)$, we shall compute $\bar L^k(\rho_z)^3$ mod $O(\eps^2)$. Thus, we have
\beq
\rho_z^3=8(8\bar z^3+12\eps A\bar z^2(\bar z+z))+O(\eps^2)
=
2^5 (2\bar z^3+3\eps A(\bar z^3+z\bar z^2))+O(\eps^2).
\eeq
We obtain  for $w=0$,
\beq
\frac1{2^k}
\bar L^k \rho_z^3
=
\left(
	2z \frac\d{\d\bar w}
	+ \epsilon \left(
		- B\bar w \frac\d{\d\bar z}
		+ A(z + \bar z) \frac\d{\d\bar w}
	\right)
\right)^k
\rho_z^3
\eeq
and since $\rho_z$ is independent of $w$,
\beq
	\bar L\rho_z^3 = \bar L^3\rho_z^4 = \bar L^4\rho_z^3 = O(\eps^2),
\eeq
and
\beq
	\frac14 \bar L^2 \rho_z^3
	= -2z \epsilon B \frac\d{\d \bar z} \rho_z^3
	= -2^7\cdot 3\epsilon B z \bar z^2 + O(\eps^2)
\eeq
%
%
%

\subsection{Terms of the form $\bar L^k\rho_z^2\rho_w$; second row.}
As before, replacing occurences of $w$ (but not $\bar w$) with $0$ (written $\cong$),
we obtain
\beq
\begin{aligned}
\rho_z^2\rho_w &\cong
	(4\bar z+2\eps A(\bar z+z))^2(4 + 2\eps B) \bar w\\
	&=2^5(2\bar z^2\bar w+\eps((2A+B)\bar z^2\bar w+2Az\bar z\bar w)+O(\eps^2),
\end{aligned}
\eeq
We compute, mod $O(\eps^2)$
\begin{multline}
\frac1{2^{k+5}}
\bar L^{k} (\rho_z^2 \rho_w) \\
\cong \left(
	2z \frac\d{\d\bar w}
	+ \epsilon \left(
		- B\bar w \frac\d{\d\bar z}
		+ A(z + \bar z) \frac\d{\d\bar w}
	\right)
\right)^k
\left(
(2\bar z^2\bar w+\eps((2A+B)\bar z^2\bar w+2Az\bar z\bar w)
\right),
\end{multline}
and hence for $w=0$,
\beq
	\bar L^{k} (\rho_z^2 \rho_w) = O(\epsilon),
	\  k\ne 1;\quad \bar L^{4} (\rho_z^2 \rho_w) = O(\epsilon^2),
\eeq
and
\beq
	\bar L (\rho_z^2 \rho_w) = 2^8 z\bar z^2 + O(\epsilon).
\eeq
%
%
%
%
%
%

\subsection{Terms of the form $\bar L^k\rho_z\rho_w^2$; third row.} We can obtain the formulas in this case by considering the previous subsection and interchanging the roles of $z$ and $w$. We obtain
\beq
\begin{aligned}
\rho_z\rho_w^2 &\cong
	(4\bar z+2\eps A(\bar z+z))(4 + 2\eps B)^2 \bar w^2\\
	&=2^5
		\left(
			2\bar z
		+\eps((2B+A)\bar z+Az)
		\right)\bar w^2
					+O(\eps^2),
\end{aligned}
\eeq
As before we obtain mod $O(\eps^2)$,
\begin{multline}
\frac1{2^{k+5}}
\bar L^{k} (\rho_z \rho_w^2) \\
\cong \left(
	2z \frac\d{\d\bar w}
	+ \epsilon \left(
		- B\bar w \frac\d{\d\bar z}
		+ A(z + \bar z) \frac\d{\d\bar w}
	\right)
\right)^k
2^5\left(
(2\bar z
		+\eps((2B+A)\bar z+Az)\right)\bar w^2,
\end{multline}
from where as before, for $w=0$,
\beq
\bar L^k ( \rho_z\rho_w^2 ) = O(\epsilon), \quad k\ne 2,
\eeq
\beq
\bar L^2 ( \rho_z\rho_w^2 ) = 2^{11}z^2\bar z +  O(\epsilon),
\eeq
and
\beq
\bar L^4 ( \rho_z\rho_w^2 )
 = 2^7
{4 \choose {2}}
\epsilon
\left(
	2z \frac\d{\d\bar w}
	- \epsilon  B\bar w \frac\d{\d\bar z}
	\right)^2
	(4\bar z z^2)
 +  O(\epsilon^2)
=  2^{11}
{4 \choose {2}}
\epsilon
B z^3
 +  O(\epsilon^2).
\eeq

%
%
%
%
%
%
%

\subsection{Terms of the form $\bar L^k\rho_w^3$; fourth row.} Following the same strategy as above, 
we obtain
\beq
\rho_w^3 \cong
(4 + 2\eps B)^3 \bar w^3 =
	2^3 (8 + 12\eps B) \bar w^3
	+O(\eps^2),
\eeq
and for $w=0$,
\beq
\bar L^k ( \rho_w^3 ) = O(\epsilon^2), \quad k\ne 3,
\eeq
\beq
\bar L^3 (\rho_w^3 ) = 2^{9} \cdot 6 z^3 +  O(\epsilon).
\eeq

%
%
%
%

\subsection{Calculation of $\eps^2$-term of $A_3(\rho_\eps)$ along $w=0$}
From our calculations in the subsections above we obtain for $w=0$:
\beq
A_3(z,\bar z)=
\begin{pmatrix}
	2^6\bar z^3+O(\epsilon) &0& O(\epsilon) &0& 	0\\
	O(\epsilon)& 8z\bar z^2 +  	O(\epsilon) & 	O(\epsilon)&	O(\epsilon)& 0 \\
	O(\epsilon)&	O(\epsilon)& 2^{11} z^2\bar z & O(\epsilon)&	2^{11}
{4 \choose {2}} \epsilon B z^3\\
	0&	0& O(\epsilon)&	2^{10}\cdot 3 z^3 & 	0\\
	2^5 \epsilon B\bar z^2 &	0&	 2^{10} \epsilon Az^2 &	0&	0
\end{pmatrix}
+ O(\epsilon^2).
\eeq
Since the last row as well as the last column each has a factor $\eps$,
we conclude
\beq\Label{det-a3}
\det A_3(\rho_\eps)|_{w=0}=\eps^2\Delta_2 +O(\eps^3),
\quad
\Delta_2 =  N AB z^9\bar z^5,
\eeq
%
where $N$ is a large positive integer.

\subsection{Umbilical points on the ellipsoids $E_\eps$.} Let $S_\eps$ be the ellipse (in the $z$-plane) obtained by intersecting $E_\eps$ with the complex line $w=0$. If both $A, B>0$, we easily conclude that the winding number $W_{S_0}(\Delta_2)$ of $\Delta_2(z,\bar z)$ around the circle $S_0$, traversed in the positive direction, equals $4$; recall that the winding number is defined by \eqref{W(R)} from which $W_{S_0}(\Delta_2)=4$ follows immediately.

\begin{proof}[Proof of Theorem $\ref{Prop-ell}$]
Let us first assume that both $A,B>0$. Since, in this case, $\Delta_2$ does not vanish on $S_0$, it follows that $Q_\eps|_{w=0}$ does not vanish on $S_\eps$ for $\eps>0$ sufficiently small, where
$Q_\eps=\det A_3(\rho_\eps)$.
 It is also clear by continuity that, for sufficiently small $\eps>0$, the winding number of $Q_\eps$ around $S_0$ coincides with that of $\Delta_2$ around $S_0$, and also around $S_\eps$. We conclude that
\beq\Label{W(Qeps)}
W_{S_\eps}(Q_\eps)=4,
\eeq
for sufficiently small $\eps>0$. Now, either the set of umbilical points $\mU\subset E_\eps$ contains points of dimension at least 2, or there is a surface $\Sigma^\eps$ in $E_\eps$ that is bounded by $S_\eps$ and meets the subset of 1-regular umbilical points $\mU_1$ transversally; indeed, we can always find even a simply connected $\Sigma^\eps$ in $E_\eps$ with $\partial\Sigma^\eps=S_\eps$, and if $\mU$ has only components of dimension 0 and 1, then small local deformations of $\Sigma^\eps$ along the intersection will result in only transversal intersections along $\mU_1$. It now follows from \eqref{W(Qeps)} and Proposition \ref{Wvsindex} that
\beq
\sum_{p\in \Sigma^\eps\cap \mU_1}\iota_{\Sigma^\eps_p}(p)=-2.
\eeq
In particular, either $\mU$ has points of dimension at least 2 or contains at least one curve of stable umbilical points when both $A,B>0$.

In the remaining case where, say, $B=0$, the ellipsoid is invariant under the circle action $(z,w)\mapsto (z,e^{it}w)$ and therefore has umbilical points along the curve of fixed points $(z,0)$, in view of the special Chern--Moser normalization at non-umbilical points \cite{CM74}, pp.\, 246--247.
\end{proof}

\section{Umbilical points on perturbations of the sphere}\Label{PertSec}

We shall consider perturbations $M_\eps\subset \bC^2$ of the unit sphere
given by $\rho=\rho^\epsilon=0$, where
\beq\Label{pert}
\rho^\eps:= \rho^0 +  \eps \rho',
\quad \rho^0 := -1 + z\bar z + w\bar w,
\eeq
$\rho'$ is a smooth real-valued function, and $\eps$ is a small real parameter.
For $\eps=0$ we recover the unit sphere $S^3=M_0$ and hence $\det A_3=0$.
For $\eps\neq 0$, we shall consider the power series expansion of $\det A_3$ in $\eps$. In that expansion, we shall compute the linear term in $\eps$.
Since the expansion of $\rho_{Z^2}$ begins with a linear term in $\eps$,
the only nonzero contribution to the linear term in $\eps$ in the determinant \eqref{a3} will come
from $0$th order terms (in $\eps$) in the first $4$ rows and $1$st order terms in the last row.
Furthermore, only $0$th order terms in the expansion of $L$ will contribute.
Thus for our computation, we only need use the terms with
$$(\rho^0_w, \rho^0_z)=(\bar w, \bar z), \quad L_0 = -\bar w \d_z + \bar z \d_w,$$
and hence the desired coefficient  of $\eps$ is
\begin{equation}\label{epscoeff}
\det
\begin{pmatrix}
D^0_3 &0\cr
* & \bar L_0^4 (\rho'_{Z^2}(L_0,L_0))
\end{pmatrix}
= (\det D^0_3) \bar L_0^4 (\rho'_{Z^2}(L_0,L_0)),
\end{equation}
where $D^0_3$ is calculated using $\rho^0$. By Proposition \ref{dn-nonvanish}, we conclude:

\bp\Label{Lin-eps} For a perturbation of the form \eqref{pert},
\beq\Label{Lin-epsterm}
\det A_3(\rho^\eps)=c_0\bar L_0^4 (\rho'_{Z^2}(L_0,L_0))\eps+O(\eps^2),
\eeq
where $c_0$ is a universal polynomial that does not vanish on the unit sphere $\rho^0=0$.
\ep
We note that
\begin{equation}\Label{rhoZ2}
\rho'_{Z^2}(L_0,L_0)= (-\bar w)^2 \rho'_{z^2} - 2\bar z\bar w \rho'_{zw} + \bar z^2 \rho'_{w^2},
 \end{equation}
and observe that the coefficients in $\bar L_0$ are holomorphic, and hence repeated applications of $\bar L_0$ will not result in any differentiations of the coefficients, and we obtain
\beq\Label{l4}
 \bar L^4_0=(-w\partial_{\bar z}+z\partial_{\bar w})^4=w^4\partial^4_{\bar z}-4zw^3\partial^3_{\bar z}\partial_{\bar w}+6z^2w^2\partial^2_{\bar z}\partial^2_{\bar w}-4z^3w\partial_{\bar z}\partial^3_{\bar w}+z^4\partial^4_{\bar w}.
\eeq
We shall consider polynomial perturbations of the form $\rho'=\sum_{k=2}^m\rho'_{k}$, where $\rho'_k$ are homogeneous polynomials of degree $k$ in $Z=(z,w)$ and $\bar Z$. We may decompose $\rho'_k$ further into bidegree, $\rho'_k=\sum_{p+q=k}\rho'_{p,q}$, where each $\rho_{p,q}$ is of bidegree $(p,q)$. Since our perturbations $\rho'$ are real-valued, we must have $\rho'_{q,p}=\overline{\rho'_{p,q}}$.

We shall use the notation $\mathcal H_k$ for the space of homogeneous polynomials of degree $k$, and $\mathcal H_{p,q}$ for those of bidegree $(p,q)$. We note that if $R\in \mathcal H_{p,q}$, then $\bar L_0^4 (R_{Z^2}(L_0,L_0))\in \mathcal H_{p+2,q-2}$. We also note that in this case $R_{Z^2}(L_0,L_0)\in \mathcal H_{p-2,q+2}$, and we conclude that  $\bar L_0^4 (R_{Z^2}(L_0,L_0))=0$ unless both $p$ and $q$ satisfy $p,q\geq 2$. Let us for brevity use the notation
\beq\Label{Q0}
Q^0(R):=\bar L_0^4 (R_{Z^2}(L_0,L_0)),
\eeq
so that
\beq\Label{Q}
Q=Q(\rho^\eps):=\det A_3(\rho^\eps)=c_0Q^0(\rho')\eps+O(\eps^2).
\eeq
We may then summarize the discussion above as follows.

\bp\Label{Q0decomp} For a real-valued polynomial $\rho'$ of degree $m$, decomposed into homogeneous components $\rho'_k$ and further decomposed into bidegree $\rho'_{p,q}$, of the form
\beq\Label{rho'decomp}
\rho'=\sum_{k=2}^m\rho'_k=\sum_{k=2}^m\sum_{p+q=k}\rho'_{p,q},\quad \rho'_{p,q}=\overline{\rho'_{q,p}},
\eeq
it holds that
\beq\Label{Q0decomp}
Q^0(\rho')=\sum_{k=4}^m\sum_{l=4}^k Q^0_{l,k-l},\quad Q^0_{l,k-l}=Q^0(\rho'_{l-2,k-l+2}).
\eeq
\ep

We have the following technical result.

\bp\Label{gengenpert}
Let $\rho'$ be a real-valued polynomial of degree $m$, and decompose $Q^0(\rho')$ as in \eqref{Q0decomp}. Assume that:
\begin{itemize}
\item[(i)] The real-algebraic variety $\mathcal V:=\{Q^0(\rho')=0\}\cap S^3$ in $S^3$ has no points of dimension $\geq 2$.
\item[(ii)] $Q^0_{l,k-l}=0$ for
$4\leq l\leq k/2$.
\end{itemize}
Then, for sufficiently small $\eps>0$, the set of umbilical points $\mU$ on the perturbation $M_\eps$ contains either points of dimension $\geq 2$ or a curve of stable umbilical points.
\ep

\begin{Rem}\Label{gengenrem} We make a few observations:
\begin{itemize}
\item If the degree $m\leq 7$, then condition (ii) is vacuous, and hence only (i) is required in this case.
\item If the degree $m\leq 3$, then condition (i) is never satisfied, since $Q^0(\rho')$ vanishes completely. In particular, as noted in the previous section, for real ellipsoid perturbations $E_\eps$ we have $Q^0(\rho')=0$. Nevertheless, as is proved in Theorem \ref{Prop-ell}, real ellipsoids do have umbilical points.
\end{itemize}
\end{Rem}

To prove Proposition \ref{gengenpert}, we need the following lemma:

\bl\Label{Excircle}
Let $\rho'$ be a real-valued polynomial such that condition {\rm (i)} in Proposition $\ref{gengenpert}$ holds. Then, there is point $Z_0=(z_0,w_0)\in S^3$ such that $Q^0(\rho')$ does not vanish on the circle $S_0\colon t\mapsto e^{it}Z_0$ in $S^3$.
\el

\begin{proof} For $Z_1:=(z_1,w_1)\in S^3$, consider the circle $S_1$ in $S^3$ parametrized by $t\mapsto e^{it}(z_1,w_1)$. Let $\Sigma\subset S^3$ be a germ at $Z_1$ of an open, real-analytic surface, transverse to $S_1$ at this point. Consider the real-analytic map $\Gamma\colon \Sigma\times S^1\to S^3$, given by $\Gamma(z,w,t)=e^{it}(z,w)$ in local coordinates $t\to e^{it}$ on $S^1$. This map realizes an open subset $\Omega$ of $S^3$ as an $S^1$-fibration over $\Sigma$. If we let $\pi\colon \Omega\to \Sigma$ be the projection, then we can consider $\pi(\mathcal V)\subset \Sigma$, where $\mathcal V$ is the zero locus of $Q^0(\rho')$ as in Proposition \ref{gengenpert}. Since $\mathcal V$, by condition (i), has no points of dimension 2 or 3, $\pi(\mathcal V)$ is a proper sub-analytic subset of the open surface $\Sigma$. Thus, by choosing $Z_0=(z_0,w_0)$ in $\Sigma$ outside this projection, we find the desired oriented circle $S_0$, parametrized by $t\to \Gamma(z_0,w_0,t)$.
\end{proof}

\begin{Rem} We may parametrize all great circles on $S^3$ by blowing up the origin in $\bC^2$. In this way, $S^3$ becomes the unit circle in the universal line bundle $O(-1)$ over $\bP^2$. The corresponding projection $\pi \colon O(-1)\to \bP^2$ is algebraic and then the set $\pi(\mathcal V)$ of unit circles to avoid is a closed semialgebraic set in $\bP^2$.
\end{Rem}

We now proceed with the  proof of Proposition \ref{gengenpert}. Let $S_0\colon t\mapsto e^{it}Z_0$, with $Z_0=(z_0,w_0)\in S^3$, be the circle provided by Lemma \ref{Excircle}, and define a polynomial $P(\zeta,\bar \zeta)$  of degree $m$ in the variable $\zeta\in\bC$ by
\beq\Label{PfromQ}
P(\zeta,\bar \zeta):=Q^0(\rho')(\zeta Z_0,\overline{\zeta Z_0}).
\eeq
By construction of $S_0$, $P$ does not vanish on the unit circle. The decomposition of $Q^0(\rho')$ into bidegree, given by \eqref{Q0decomp}, yields a decomposition of $P$ into bidegree:
\beq
P(\zeta,\bar\zeta)=\sum_{k=4}^m\sum_{l=4}^kp_{l,k-l}
\zeta^l\bar \zeta^{k-l},\quad p_{l,k-l}
:=Q_{l,k-l}(\zeta Z_0,\overline{\zeta Z_0})/\zeta^l\bar \zeta^{k-l}.
\eeq
On the unit circle $\zeta=e^{it}$, $P(\zeta,\bar \zeta)$ coincides with a rational function $R(\zeta)$,
\beq
R(\zeta)=P(\zeta,1/\zeta)=\frac{\sum_{k=4}^m\sum_{l=4}^k
p_{l,k-l}\zeta^{2l+m-k}}{\zeta^m}.
\eeq
If we define
\beq
p(\zeta):=\sum_{k=4}^m\sum_{l=4}^k
p_{l,k-l}\zeta^{2l+m-k},
\eeq
then by the construction of $p(\zeta)$ and the argument principle we conclude:

\bl\Label{ArgP}
Let $n$ denote the number of zeros (counted with multiplicities) of $p(\zeta)$ in the unit disk $|\zeta|<1$. Then
\beq
W_{S_0}(Q^0)=n-m,
\eeq
where $Q^0=Q^0(\rho')$ and $S_0$ is the circle in the construction of $p(\zeta)$ above.
\el

We may now complete the proof of Proposition \ref{gengenpert}.

\begin{proof}[Proof of Proposition $\ref{gengenpert}$] Recall that the set of umbilical points $\mU$ on $M_\eps$ is given by $Q=0$, where $Q$ is as in \eqref{Q}. Let $S_0$ the circle on $S^3$ as above, and let $S_\eps$ be the perturbed oriented curve on $M_\eps$ obtained as the intersection between the complex subspace through $Z_0=(z_0,w_0)\in S^3$ and $M_\eps$. As in the proof of Theorem \ref{Prop-ell}, for sufficiently small $\eps>0$, we have
\beq
W_{S_\eps}(Q)=W_{S_0}(Q^0)=n-m,
\eeq
where $n$ and $m$ are as in Lemma \ref{ArgP}. We claim that $n-m\neq 0$. Indeed, by condition (ii), the coefficients $p_{l,k-l}=0$ for $l\leq k/2$. If we let $r$ denote the minimum integer $r=2l+m-k$ for which $p_{l,k-l}\neq 0$, then $p(\zeta)$ is divisible by $\zeta^{r}$ and since $r>m$, we conclude that $n>m$. The proof of Proposition \ref{gengenpert} is now completed in the same way as the proof of Theorem \ref{Prop-ell}.
\end{proof}

\subsection{The sphere $S^3$ as a circle bundle over $\bP^1$} We recall here the idea of realizing the sphere as the unit circle in the universal bundle
$\pi\colon J:=O(-1)\to \bP^1$; this idea has been extended to more general three-dimensional CR manifolds with a CR circle action by Bland--Duchamp \cite{BlandDuchamp91} and Epstein \cite{Epstein92}. Recall that $J$
naturally embeds into $\bP^1\times \bC^2$ in such a way that the fiber $J_Z$ over a point $Z$ in homogeneous coordinates, $Z=[z\colon w]\in \bP^1$, is the complex line through $(z,w)\in \bC^2$; $J_Z$ is parametrized by $\zeta\mapsto \zeta(z,w)$. The standard metric $|\cdot|$ on $J$ is the one induced by the Euclidian metric on $\bC^2$; if $s(Z)$ is a non-vanishing local section in $J$ and we write $s(Z)=(u,v)\in \bC^2$, then $|s|^2:=|u|^2+|v|^2$. The unit circle bundle $\tilde S^3:=\{\lambda\in J\colon |\lambda|^2=1\}$ is CR isomorphic to the unit sphere $S^3$ in $\bC^2$. Indeed, if we view the total space $J$ as the blow-up of the origin in $\bC^2$, then the CR isomorphism $\tilde \pi|_{\tilde S^3}\colon \tilde S^3\to S^3$ is the blow-down map $\tilde \pi\colon J\to \bC^2$ restricted to $\tilde S^3$. For convenience, we shall simply identify $S^3$ with $\tilde S^2$ via this isomorphism; in this identification, the fibers $\pi^{-1}(Z)$ in $\tilde S^3$ correspond to the great circles $t\mapsto e^{it}Z$.

We note that if $\rho'$ is a real-valued polynomial, then the projection $\pi(\mathcal V)\subset \bP^1$ of the real-algebraic subvariety $\mathcal V\subset S^3\cong \tilde S^3$, defined to be the zero locus of $Q^0=Q^0(\rho')$ as in Proposition \ref{gengenpert}, is a closed semialgebraic subset. An inspection of the proof of Proposition \ref{gengenpert} reveals immediately that condition (i) in the assumptions of this proposition can be replaced by the assumption that $\pi(\mathcal V)\neq \bP^1$. As is shown in Lemma \ref{Excircle}, condition (i) implies $\pi(\mathcal V)\neq \bP^1$, and the latter property is the only one used in the proof of Proposition \ref{gengenpert}. For convenience, we state the result here.

\bp\Label{gengenpert'}
Let $\rho'$ be a real-valued polynomial of degree $m$, and decompose $Q^0(\rho')$ as in \eqref{Q0decomp}. Assume that:
\begin{itemize}
\item[(i$'$)] $\bP^1\setminus \pi(\mathcal V)\neq \emptyset$, where $\mathcal V\subset S^3\cong\tilde S^3$ and $\pi\colon \tilde S^3\to \bP^1$ are as above.
\item[(ii)] $Q^0_{l,k-l}=0$ for $4\leq l\leq k/2$.
\end{itemize}
Then, for sufficiently small $\eps>0$, the set of umbilical points $\mU$ on the perturbation $M_\eps$ contains either points of dimension $\geq 2$ or a curve of stable umbilical points.
\ep

\subsection{Generic perturbations of the sphere} We shall denote by $\mathcal P_m$ the space of all polynomials in $Z=(z,w)$ and $\bar Z$ of degree at most $m$, and by $\RP_m$ the real subspace of those that are real-valued. Thus, we have $$\mathcal P_m=\bigoplus_{p+q\leq m}\mathcal H_{p,q}$$ and $\rho'\in\mathcal P_m$ belongs to $\RP_m$ when $\rho'_{p,q}=\overline{\rho_{q,p}}$ for all $p,q$. We shall show that condition (i$'$) in Proposition \ref{gengenpert'} is generic. More precisely, we shall prove the following:

\bp \Label{gencondi} The set $\Pi_m$ of polynomials $\rho'$ in $\RP_m$ such that $\pi(\mathcal V)=\bP^1$, where $\mathcal V\subset S^3\cong\tilde S^3$ and $\pi\colon \tilde S^3\to \bP^1$ are as in Proposition $\ref{gengenpert'}$, is a real-analytic subvariety in $\RP_m$. Moreover, if $\mathcal A\subset \RP_m$ is any real subspace containing $\mathcal A_p:=\{ez^p\bar w^p+\bar e w^p\bar z^p\colon e\in \bC\}$ for some $2\leq p\leq m/2$, then $\Pi_m\cap \mathcal A$ has strictly smaller dimension than $\mathcal A$.
\ep

\begin{proof} Let $\tilde z=z/w$ be a local coordinate in the chart $U_0=\{[z\colon w]\in \bP^1\colon w\neq 0\}$ in $\bP^1$ and $\tilde Q^0=\tilde Q^0(\tilde z,\bar{\tilde z};\zeta,\bar\zeta)$ the polynomial $Q^0=Q^0(\rho')$ for some $\rho'\in \RP_m$ in the local trivialization
$$U_0\times\bC=\bC\times\bC\cong J|_{U_0}\subset U_0\times\bC^2,$$
given by
$$
(\tilde z,\zeta)\mapsto (\tilde z;\tilde\pi(\tilde z,\zeta)),\ \tilde\pi(\tilde z,\zeta):=\zeta(\tilde z,1).
$$
In other words, $\tilde Q^0=Q^0\circ \tilde \pi$;  we shall denote by $\tilde Q^0_{p,q}=Q^0_{p,q}\circ\tilde \pi$, so that we have the decomposition (see Proposition \ref{Q0decomp})
\beq
\tilde Q^0=\sum_{k=4}^m\sum_{l=4}^k\tilde Q^0_{l,k-l},
\eeq
where each component $\tilde Q^0_{p,q}$ takes the form
\beq
\tilde Q^0_{p,q}(\tilde z,\bar{\tilde z};\zeta,\bar\zeta)=q_{p,q}(\tilde z,\bar{\tilde z})\zeta^p\bar \zeta^q=q_{p,q}(\tilde z,\bar{\tilde z})\zeta^{p-q}|\zeta|^{2q},
\eeq
with
\beq
q_{p,q}(\tilde z,\bar{\tilde z})=\tilde Q^0_{p,q}((\tilde z,1),(\bar{\tilde z},1))=
\left(\sum_{\alpha\leq p,\, \gamma\leq q}c_{pq;\alpha\bar \gamma}\tilde z^{\alpha}\bar{\tilde z}^{\gamma}\right).
\eeq
for suitable coefficients $c_{pq;\alpha\bar\gamma}$.
Recall that $\tilde S^3\subset J$ is given in these coordinates by
\beq\Label{tildeS^3}
|\zeta|^2(1+|\tilde z|^2)=1.
\eeq
Consequently, each $\tilde Q^0_{p,q}$ coincides on $\tilde S^3$ with the function
\beq
R_{p,q}(\tilde z,\bar{\tilde z};\zeta,\bar\zeta)=\frac{q_{p,q}(\tilde z,\bar{\tilde z})}{(1+|\tilde z|^2)^q}\,\zeta^{p-q},
\eeq
and $\tilde Q^0$ coincides with $R$, where
\beq\Label{R1}
R=\sum_{k=4}^m\sum_{l=4}^kR_{l,k-l}=\sum_{k=4}^m\sum_{l=4}^k\frac{q_{l,k-l}(\tilde z,\bar{\tilde z})}{(1+|\tilde z|^2)^{k-l}}\,\zeta^{2l-k},
\eeq
a rational function in $\zeta$ with coefficients that are rational functions in $\tilde z$ and $\bar{\tilde z}$. Note that the powers of $\zeta$ range from $8-m$ to $m$. Let us collect terms of equal powers in $\zeta$ and rewrite $R$ in \eqref{R1} in the form
\beq\Label{R2}
R=\sum_{r=8-m}^m \frac{b_r(\tilde z,\bar{\tilde z})}{(1+|\tilde z|^2)^{s_r}}\,\zeta^r
=\frac{1}{\zeta^{m-8}}\sum_{r=0}^{2m-8} \frac{b_{r+8-m}(\tilde z,\bar{\tilde z})}{(1+|\tilde z|^2)^{s_{r+8-m}}}\,\zeta^r,
\eeq
where the $s_r$ are (easily computable but not important) positive integers, and the $b_r$ are polynomials in $(\tilde z,\bar{\tilde z})$.

Now, by definition of the set $\Pi_m$, we have $\rho'\in \Pi_m$ precisely when $R$ as a rational function in $\zeta$ has at least one root on the circle \eqref{tildeS^3} for every $\tilde z\in U_0\subset \bP^1$. Observe that that set $B_k$ of coefficients $a=(a_0,\ldots,a_k)\in \bC^{k+1}$ such that the polynomial $a_0+a_1\zeta+\ldots+a_k\zeta^k$ has a root on the unit circle forms a real-algebraic, Levi flat (singular) hypersurface. Thus, $\rho'\in \Pi_m$ translates into the condition that
\beq\Label{br1}
\left(\frac{b_{8-m}(\tilde z,\bar{\tilde z})}{(1+|\tilde z|^2)^{s_{m-8}}},\frac{b_{9-m}(\tilde z,\bar{\tilde z})}{(1+|\tilde z|^2)^{s_m+1/2}},\ldots,
\frac{b_m(\tilde z,\bar{\tilde z})}{(1+|\tilde z|^2)^{s_m+(2m-8)/2}}\right )\in B_{2m-8},\quad \forall \tilde z\in U_0.
\eeq
By unraveling the construction of $R$, we note that if we expand $\rho'$ in the monomial basis $Z^I=z^\alpha w^\beta$ of $\mathcal P_m$, i.e.,
\beq
\rho'=\sum_{|I|+|J|\leq m}e_{I\bar J}Z^I\bar Z^J,\quad e_{I\bar J}=\overline{e_{J\bar I}},
\eeq
then the components in \eqref{br1}
$$
\frac{b_r(\tilde z,\bar{\tilde z})}{(1+|\tilde z|^2)^{s_r+r/2}}
$$
are linear in $e_{I\bar J}$ and $\overline{e_{I\bar J}}$. Consequently, we deduce from the above discussion and \eqref{br1} that $\Pi_m$ is a real-algebraic subvariety in $\RP_m$.

To complete the proof of Proposition \ref{gencondi}, we must show that if $\mathcal A$ is as in the statement of the proposition, then the dimension of $\Pi_m\cap \mathcal A$ is strictly less than that of $\mathcal A$. For this, it suffices to show that $\Pi_m\cap\mathcal A\neq \mathcal A$. To this end, we compute $Q^0(z^p\bar w^p)$, for $p\geq 2$,
\beq
Q^0(z^p\bar w^p)=\bar L_0^4(p(p-1)z^{p-2}\bar w^{p+2})=(p+2)(p+1)p^2(p-1)^2z^{p+2}\bar w^{p-2},
\eeq
and similarly,
\beq
Q^0(w^p\bar z^p)=\bar L_0^4(p(p-1)w^{p-2}\bar z^{p+2})=(p+2)(p+1)p^2(p-1)^2w^{p+2}\bar z^{p-2}.
\eeq
Thus, if $\rho'$ is any polynomial in $\mathcal A$, resulting in the polynomial $R$ as in \eqref{R2}, then $\rho'+ez^p\bar w^p+\bar ew^p\bar z^p$, which is also in $\mathcal A$ for all $e\in \bC$, results in
\beq
R'=R+\frac{(e\tilde z^{p+2}+\bar e\bar {\tilde z}^{p-2})}{(1+|\tilde z|^2)^{p-2}}\, \zeta^4.
\eeq
From this we easily deduce that if $\rho'\in \Pi_m$, then $\rho'+ez^p\bar w^p+\bar ew^p\bar z^p$ will not be in $\Pi_m$ for $e\neq 0$; indeed, since
\beq
\tilde z\mapsto \frac{(e\tilde z^{p+2}+\bar e\bar {\tilde z}^{p-2})}{(1+|\tilde z|^2)^{p-2}},\quad e\neq 0,
\eeq
maps onto an open neighborhood of $0$ in $\bC$, this statement follows from the following simple observation:

\bl If $p(\zeta)=\zeta^n+a_{n-1}\zeta^{n-1}+\ldots +a_0$ has a root on the unit circle, then the set of $b\in \bC$ such that $p(\zeta)+b\zeta^k$ has a root on the unit circle is a real-algebraic, possibly singular curve (real-algebraic variety of dimension one).
\el

\begin{proof}  Consider the (symmetric) finite polynomial mapping $\Phi\colon \bC^{n}\to \bC^{n}$ sending a collection of roots $\tau=(\tau_1,\ldots,\tau_{n})$ to the collection of coefficients $a=(a_0,\ldots, a_{n-1})$ of the polynomial
\beq
p(\zeta)=\zeta^n+a_{n-1}\zeta^{n-1}+\ldots +a_0:=(\zeta-\tau_1)\ldots(\zeta-\tau_n).
\eeq
Pick $p_0(\zeta)$ such that one its roots is on the unit circle, i.e., $a^0=\Phi(\tau^0)$ with $\tau^0$ in the Levi flat (singular) hypersurface $H=\cup_{j=1}^n H_j$, with $H_j:=\{\tau\colon |\tau_j|=1\}$. The polynomials $p_e(\zeta):=p_0(\zeta)+b\zeta^k$ correspond to points $a^b=a^0+(0,\ldots,b,\ldots,0)$ (with $b$ in the $(k+1)$th component) and hence their roots $\tau^b$ belong to the complex 1-dimensional subvariety $\Phi^{-1}(X_k)$, where $X_k$ denotes the complex curve $b\mapsto  a^0+(0,\ldots,b,\ldots,0)$. We claim that $\Phi^{-1}(X_k)$ is not contained in $H$, which will prove the conclusion of the lemma. Indeed, $\Phi^{-1}(X_k)$ could only be contained in the Levi flat $H$ if it were contained in one of its leaves $\tau_j=c$, with $c$ constant, which is clearly impossible.
\end{proof}
As mentioned above, we have now shown that the real-algebraic subvariety $\Pi_m$ satisfies $\Pi_m\cap\mathcal A\neq \mathcal A$, which completes the proof of Proposition \ref{gencondi}.
\end{proof}

\subsection{Generic perturbations of almost circular type} Recall that a real hypersurface $M\subset \bC^{2}$ is called {\it circular} if $Z\in M$ implies $e^{it}Z\in M$ for all $e^{it}\in S^1$. For perturbations $M_\eps$ of the sphere, as in \eqref{pert}, it is straightforward to verify that the $M_\eps$ are circular for all sufficiently small $\eps>0$ if and only if in the decomposition \eqref{rho'decomp} we have $\rho'_{p,q}=0$ for $|p-q|\neq 0$. It was shown in \cite{EDumb15} that compact, circular real hypersurfaces in $\bC^2$ always have umbilical points. Here we shall consider perturbations $M_\eps$ that are {\it almost circular}, which we define to be those for which, in the decomposition \eqref{rho'decomp} of $\rho'$, we have $\rho'_{p,q}=0$  when $|p-q|\geq 4$; we also say that such $\rho'$ are almost circular. We easily observe that a polynomial $P=P(Z,\bar Z)$ is almost circular if and only if its Fourier coefficients $\hat P_k$ vanish for $|k|\geq 4$:
$$
\hat P_{k}(Z,\bar Z):=\frac{1}{2\pi}\int_{0}^{2\pi}P(e^{it}Z,e^{-it}\bar Z)e^{-ikt}dt=0,\quad |k|\geq 4.
$$

Recall that $\mathcal P_m$ denotes the space of all polynomials in $Z=(z,w)$ and $\bar Z$ of degree at most $m$. We shall denote by $\AC_m$ the real subspace of those that are real-valued and almost circular. Thus, $\rho'\in\mathcal P_m$ belongs to $\AC_m$ when $\rho'$ is real-valued (i.e., $\rho'\in \RP_m$) and $\rho'_{p,q}=0$ for $|p-q|\geq 4$. We note that $\mathcal A=\AC_m$ satisfies the hypothesis in Proposition \ref{gencondi} for all $m\geq 2$ and with any $2\leq p\leq m/2$.

\bt\Label{ACpertThm}
For $m\geq 4$, there is a real-algebraic subvariety $\Xi_m\subset \AC_m$ of dimension strictly less than that of $\AC_m$ such that if $\rho'\in \AC_m\setminus\Xi_m$, then, for sufficiently small $\eps>0$, the set of umbilical points $\mU$ on the perturbation $M_\eps$, given by \eqref{pert}, contains either points of dimension $\geq 2$ or a curve of stable umbilical points.
\et

\begin{proof} We shall let $\Xi_m$ be $\Xi_m:=\Pi_m\cap \AC_m$, where $\Pi_m$ is as  defined in Proposition \ref{gencondi}. As noted above, $\mathcal A=\AC_m$ satisfies the hypotheses in Proposition \ref{gencondi} and, hence, we conclude that $\Xi_m$ is a real-algebraic subvariety of strictly lower dimension that $\AC_m$. The conclusion of Theorem \ref{ACpertThm} now follows from Proposition \ref{gengenpert'}, since $\rho'\in \AC_m$ clearly guarantees that condition (ii) in that proposition holds; indeed, for $\rho'_{p,q}$, we have $Q^0(\rho'_{p,q})=Q^0_{p+2,q-2}$ and if $|p-q|\leq 3$, then $l=p+2\geq(p+q+1)/2>k/2$, which is the requirement in condition (ii).
\end{proof}

\begin{Rem}\Label{ACrem}
\begin{itemize}
\item
Recall that if, for example, $m=2p$ and
$$
\rho'=\rho'_{p-1,p+1}+\rho_{p,p}+\rho'_{p+1,p-1},\quad  \rho'_{p+1,p-1}=\overline{\rho'_{p-1,p+1}}
$$
($\implies \rho'\in \AC_m$), then
$$
Q^0=Q^0(\rho')=Q^0_{p+1,p-1}+Q^0_{p+2,p-2}+Q^0_{p+3,p-3}.
$$
We note that there are plenty of polynomials of this form,
$$
Q=Q_{p+1,p-1}+Q_{p+2,p-2}+Q_{p+3,p-3},
$$
such that $\pi(\mathcal V)=\bP^1$, where $\mathcal V$ denotes the zero locus of $Q$ in $\tilde S^3$. For example, any $Q$ of the form
$$
Q=(z+\bar z)(Q'_{p-1,p}+Q'_{p,p-1})
$$
will satisfy this, as the reader can easily verify. However, we do not know any non-trivial examples of such $Q$ that are also in the image of the linear map $Q^0$, i.e., of the form $Q=Q^0(\rho')$ with $\rho'\in\AC_m$.
\item It is clear from the calculations in Section \ref{EllSec} that the real ellipsoids $E_\eps$ are not generic in the sense of Theorem \ref{ACpertThm}, i.e., these belong to $\Pi_m$.
\end{itemize}
\end{Rem}


\def\cprime{$'$}

\end{document}